\documentclass[12pt]{article}
\usepackage{amsmath, amsthm, amscd, amsfonts, amssymb, graphicx, color}
\usepackage[bookmarksnumbered, colorlinks, plainpages]{hyperref}
\newtheorem{theorem}{Theorem}[section]
\newtheorem{lemma}[theorem]{Lemma}
\newtheorem{proposition}[theorem]{Proposition}
\newtheorem{corollary}{Corollary}[theorem]
\theoremstyle{definition}
\newtheorem{definition}[theorem]{Definition}

\theoremstyle{remark}
\newtheorem{remark}{Remark}[section]
\numberwithin{equation}{section}

\begin{document}
	
	\setcounter{page}{1}
	
	%%%%%%%%%%%%%%%%%%%%%%%%%%%%
	
	\begin{center}
		{\Large \textbf{On Some Characterization of $GS$-exponential kind of Convex Functions}}
		
		\bigskip

		\textbf{Ehtesham Akhter$^a$} and \textbf{Musavvir Ali$^{b,*}$}\\
		\textbf{}  \textbf{}\\

		{\small $^{a,~b}$  Department of Mathematics,\\ Aligarh Muslim University, Aligarh-202002, India}
		
		\noindent
		{
			E-mail addresses: $b^*$musavvir.alig@gmail.com (Corresponding author),\\ ehteshamakhter111@gmail.com.}
	\end{center}
	\bigskip

	{\abstract
		\noindent
		This manuscript introduces the idea of $GS$-exponential kind of convex functions and some of their algebraic features, and we introduce a new class $GS$-exponential kind of convex sets. In addition, we describe certain fundamental  $GS$-exponential kind of convex function with characteristics in both the general and the differentiable cases. We establish the sufficient conditions of optimality and offer the proof for unconstrained as well as inequality-constrained programming while considering the assumption of $GS$-exponential kind of convexity.}\\
	
	
	
	\section{Introduction}
	Due to the importance of convexity and its generalisations in the study of optimality to resolve mathematical issues, researchers have concentrated a lot of their efforts on generalised convex functions for this purpose. As an illustration, Hudzik and Maligranda (1994) \cite{Hudzik}, investigated at two distinct forms of $s$-convexity and found that $s$-convexity in the next meaning is basically more significant than in the first sense whenever $(0< s< 1)$. Youness (1999) \cite{Youness} expanded the definitions of convex sets and functions to create a new class of sets and functions known as $E$-convex sets and $E$-convex functions.
	Yang (2001) \cite{X}  enhanced Youness's paper \cite{Youness} by incorporating certain illustrations.
	
	In recent years, academic experts have given these generalized convex functions in additional consideration. The semi-preinvex functions were studied by X.J. Long and J.W. Peng in 2006 \cite{Long} as a generalization of the semi-preinvex functions and the $b$-vex functions. Y. Syau et al. (2009)\cite{Syau}  developed the $E$-$b$-vex function family, a novel class of functions which are the generalizations of $b$-vex functions and $E$-vex functions. In 2011, T. Emam investigated a novel class of functions known as approximately $b$-invex functions. He also discussed some of its properties and discovered the necessary optimality conditions for nonlinear programming using these functions. In their investigation of a novel class of generalized sub-$b$-convex functions and sub-$b$-convex sets, M.T. Chao et al. (2012) \cite{Chao} showed the conditions for the existence of optimal solutions for both unconstrained and inequality-constrained sub-$b$-convex programming.
	
	The study in our paper aims to introduce a new class of generalized exponential kind of convex functions termed as $GS$-exponential kind of convex functions and explores certain characteristics of the same class. This paper draws inspiration from a number of research papers \cite{Butt,Fakhar,Fulga,IH,Kadakal1,Mishra,ozcan,Shi,Wang,Zhao}. Additionally, we offer the adequate $GS$-exponential kind of convexity-derived criteria of optimality for programming with variables which are both unconstrained and inequality-constrained.

	\section{Preliminaries}
	We will go through the definitions of sub-$b$-$s$-convexity, exponential kind of convexity, and $s$-convexity of functions in this section of the manuscript. For the remainder of this work, let $V$ stand for any non-empty convex subset in $\mathbb{R}^n$. 
	\begin{definition}\label{d3} \cite{Liao}
		The function $Q: V \rightarrow \mathbb{R}$ is known as sub-$b$-$s$-convex in the second sense associated with the map $G: V \times V \times (0,1]\rightarrow \mathbb{R}$, if
		$$	Q(am_1+(1-a)m_2)\leq a^sQ(m_1)+(1-a)^sQ(m_2)+G(m_1,m_2,s)$$
		holds for all $m_1,m_2\in V, a \in [0,1 ]$ and for any fixed $s \in (0,1]$.
	\end{definition} 
	\begin{definition}\label{d4}\cite{Hudzik}
		The function $Q: V \rightarrow \mathbb{R}$ is known as $s$-convex in the second sense, if  for all $m_1,m_2\in V, a \in [0,1 ]$ and for any fixed $s \in (0,1]$, we have
		$$	Q(am_1+(1-a)m_2)\leq a^sQ(m_1)+(1-a)^sQ(m_2)$$
	\end{definition}
	\begin{definition}\label{d5}\cite{Kadakal}
		A positive function $Q: V \rightarrow \mathbb{R}$ is known as exponential kind of convex function, if
		$$Q(am_1+(1-a)m_2)\leq(e^a-1)Q(m_1)+(e^{1-a}-1)Q(m_2)$$
		holds for all $m_1,m_2\in V, a \in [0,1 ]$.
	\end{definition}
	The concepts defined as in \ref{d3}, \ref{d4} and \ref{d5}, motivates us to explore a new idea known as $GS$-exponential kind of convex function.
	\section{Main Results}
	
	\begin{definition}\label{d1}
		The function $ Q :  V\rightarrow \mathbb{R}$ is known as $GS$-exponential kind of convex function on $V$ associated with the map $G : V \times V \times (0, 1] \rightarrow \mathbb{R}$, if
		\begin{equation}\label{i1}
			Q(am_1+(1-a)m_2)\leq(e^a-1)^sQ(m_1)+(e^{1-a}-1)^sQ(m_2)+aG(m_1,m_2,s)
		\end{equation} holds for each $m_1, m_2 \in V, a \in [0,1]$ and for any fixed $s \in(0,1].$
	\end{definition}
	\begin{remark}\label{r1}
		If we take $s=1$, $Q(m_1)$ is non-negative and $G(m_1,m_2,s)=0$, the $GS$-exponential kind of convex function reduces to be exponential kind of convex function.
	\end{remark}

	\begin{theorem}\label{t1}
		If $Q_1, Q_2 : V \rightarrow\mathbb{R}$ are $GS$-exponential kind of convex function associated with the map $G_1, G_2$ respectively, then $Q_1+Q_2$ and $\beta Q_1$, $(\beta \geq 0)$ are also a $GS$-exponential kind of convex function.
	\end{theorem}
	\begin{corollary}\label{c1}
		If $ Q_i : V \rightarrow \mathbb{R}$, $(i=1,2,....., n)$ are $GS$-exponential kind of convex function associated with the map $G_i : V \times V \times (0,1 ] \rightarrow \mathbb{R},$ $(i=1,2,....,n)$, respectively, then $Q=\sum_{i=1}^{n}\beta_iQ_i, \beta\geq0, (i=1,2,...,n)$ is $GS$-exponential kind of convex function associated with the map $G=\sum_{i=1}^{n} \beta_iG_i$.
	\end{corollary}
	
	\begin{lemma}\label{l1}
		For all $a \in [0,1]$ and $s \in (0,1]$, the inequalities $(e^a-1)^s \geq a$ and  $(e^{1-a}-1)^s \geq 1-a$ hold.
	\end{lemma}

	\begin{proposition}\label{p1}
		Every convex function is $GS$-exponential kind of convex function if it has a map $G$ associated with it that is non-negative.
	\end{proposition}

	\begin{theorem}\label{t2}
		If $Q: V \rightarrow \mathbb{R}$ is the GS-exponential kind of convex function associated with the map $G$ and $ S : \mathbb{R} \rightarrow \mathbb{R}$ is a non-negative function in addition to being linear, then $S \circ Q$ is a GS-exponential kind of convex function associated with the map $S \circ G$.
	\end{theorem}

	\begin{definition}\label{d2}
		Assume that $U$ be a non-empty subset of $\mathbb{R}^{n+1}$. Then, $U$ is known as $GS$-exponential kind of convex set associated with the map $G : \mathbb{R}^{n} \times \mathbb{R}^{n} \times (0,1] \rightarrow \mathbb{R}$ if for all $(m_1,\alpha_1),(m_2, \alpha_2) \in U, m_1, m_2 \in \mathbb{R}^n, a \in [0, 1]$ and some fixed $ s \in (0,1 ],$ we have  $$(am_1+(1-a)m_2,(e^a-1)^s\alpha_1+(e^{1-a}-1)^s\alpha_2+aG(m_1,m_2,s)) \in U.$$ 
	\end{definition}
	Now, we provide a characterization of $GS$-exponential kind of convex function  $ Q : V \rightarrow \mathbb{R}$ based on their respective epigraphs, given by $$E(Q)=\{(m, \alpha): m \in V, \alpha \in \mathbb{R}, Q(m)\leq \alpha\}.$$
	\begin{theorem}\label{t3}
		A function $ Q : V \rightarrow \mathbb{R}$ is a $GS$-exponential kind of convex function associated with the map $G : V \times V \times (0,1 ] \rightarrow \mathbb{R}$, if and only if $E(Q)$ is a $GS$-exponential kind of convex set associated with the map $G$.
	\end{theorem}

	\begin{theorem}\label{t4}
		Assume that $m_2>0$ and $Q_\beta: [m_1,m_2]\rightarrow \mathbb{R}$ is a family of numerical functions associated with the map $G_\beta$ and each $G_\beta$ is a $GS$-exponential kind of convex functions and each $G_\beta$ is bounded function, also assume that $Q(m)=\sup_\beta Q_\beta(m)$ and $G(m_1,m_2,s)=\sup_\beta G_\beta(m_1,m_2,s)$. If the set (non-empty) $K=\{r \in [m_1,m_2] | Q(r)<\infty\}$, then $K$ is an interval and $Q$ is $GS$-exponential kind of convex function on $K$.
	\end{theorem}

	\begin{theorem}\label{t5}
		Let  $Q:[m_1,m_2] \rightarrow \mathbb{R}$ be a $GS$-exponential kind of convex function associated with the map $G: [m_1,m_2] \times [m_1,m_2] \times (0,1] \rightarrow \mathbb{R}$ and also let $G(m_1,m_2,s)$ is bounded, then $Q$ is also bounded on $[m_1,m_2].$
	\end{theorem}

	In this section, $Q$ is considered to be a differentiable function and $ s,a \in (0,1].$
	\begin{theorem}\label{t6}
		Let $Q: V \rightarrow\mathbb{R}$ be a non-negative differentiable $GS$-exponential kind of convex function associated with the map $G$. Then
		
		$$ (i) \nabla Q(m_2)^T(m_1-m_2)< \dfrac{(e^a-1)^s}{a}Q(m_1)+\dfrac{e^{(1-a)s}}{a}Q(m_2)+G(m_1,m_2,s)-\dfrac{o(a)}{a},$$
		$$
		(ii)\nabla Q(m_2)^T(m_1-m_2)<\dfrac{(e^s-1)^s(Q(m_1)-Q(m_2))+3Q(m_2)-o(a)}{a}+G(m_1,m_2,s)
		$$
	\end{theorem}

	\begin{theorem}\label{t7}
		Let $Q: V \rightarrow \mathbb{R}$ be a non-positive differentiable $GS$-exponential kind of convex function associated with the map $G$. Then
		$$ 	\nabla Q(m_2)^T(m_1-m_2)\leq\dfrac{(e^a-1)^s}{a}[Q(m_1)-Q(m_2)]+G(m_1,m_2,s)-\dfrac{o(a)}{a}.$$
	\end{theorem}

	\begin{corollary}\label{c2}
		Assume that $Q:V \rightarrow \mathbb{R}$ is a positive differentiable $GS$-exponential kind of convex function, then
		\begin{eqnarray*}
			\nabla [Q(m_2)-Q(m_1)]^T(m_1-m_2)&<& \dfrac{(e^a-1)^s}{a}[Q(m_1)+Q(m_2)]+\dfrac{e^{(1-a)s}}{a}[Q(m_1)+Q(m_2)]\\&&+G(m_1,m_2,s)+G(m_2,m_1,s)-2\dfrac{o(a)}{a}.
		\end{eqnarray*}
		In case if $Q$ is a negative valued, then
		$$	\nabla [Q(m_2)-Q(m_1)]^T(m_1-m_2)\leq G(m_1,m_2,s)+G(m_2,m_1,s)-2\dfrac{o(a)}{a}.$$
	\end{corollary}

	The following methods are then utilized to apply the above outcomes to nonlinear programming. So, we take the unconstrained problem (S).
	\begin{equation}\label{i22}
		(S):	\min \{Q(m), m \in V\}
	\end{equation}

	\begin{theorem}\label{t8}
		Let $ Q: V \rightarrow \mathbb{R}$ be a positive differentiable $GS$-exponential kind of convex function associated with the map $G$. Also, suppose that $m \in V$ and the inequality
		\begin{equation}\label{i10}
			\nabla Q(m)^T(n-m)> G(n,m,s)+\dfrac{3Q(m)-o(a)}{a}
		\end{equation}	
		holds for each $n \in V, a \in(0,1)$, and for any particular $s \in (0,1],$, then $n$ is the solution optimal to the problem \eqref{i22} related to $Q$ on $V$.
	\end{theorem}

	The following example of unconstrained programming is taken into consideration

\end{document}